\documentclass[12pt]{article}

\usepackage[cp1251]{inputenc}

\usepackage{amsfonts}
\usepackage{amssymb}

\makeatletter
\renewcommand{\section}{\paragraph}

\makeatother
\begin{document}
\newtheorem{lemma}{Lemma}
\newtheorem{Theorem}{Theorem}

\begin{center}
{\Large On One Uniqueness Theorem for M. Rietz Potentials}

\bigskip {\large Konstantin Izyurov\footnote{The author is partially supported by the RFBR grant 06-01-00313.} }
\end{center}

\bigskip \noindent{\small We prove that there exists a nonzero
holderian function $f:\mathbb{R}\rightarrow\mathbb{R}$ vanishing
together with its M. Rietz potential $f\ast
\frac{1}{|x|^{1-\alpha}}$ in all points of some set of positive
length. This result improves the one of D. Beliaev and V. Havin
\cite{BH}}

\bigskip\noindent{\bf0. Introduction.} Let $\alpha$ be a real
number, $0<\alpha<1$. Let $f:{{\mathbb R}\rightarrow {\mathbb C}}$
be a locally summable function satisfying the following condition:
\begin{equation}
\label{condition1} \int_{ {\mathbb R}}{\frac{|f(x)|dx}{1+|x|^{1-
\alpha}}}< +\infty.
\end{equation}
Let
$$ (U_\alpha f)(t):=\int_{{{\mathbb
R}}}{\frac{f(x)dx}{|t-x|^{1- \alpha}}},\;t\in\mathbb{R}.
$$
In case (\ref{condition1}) the function $U_\alpha f$ is defined
a.e. on $\mathbb{R}$. We call it \textit{the M. Rietz Potential},
and we call $f$ \textit{the density} of this potential. We write
$dom U_\alpha$ (the domain of $U_\alpha$) for the set of all
locally summable functions satisying (\ref{condition1}).

Let $V\in\mathbb{R}$ be a measurable set; we denote its length as
$|V|$. The uniqueness theorem mentioned in the title states that
\textit{if $f$ satisfies (\ref{condition1}) and H\"{o}lder's
condition with an exponent more than $1-\alpha$ in some
neighboarhood of $V$, while $|V|>0$ and}
\begin{equation}
\label{allnull} f|_V=U_\alpha|_V=0,
\end{equation}
\textit{then $f=0$ a.e. on $\mathbb{R}$.}

This theorem follows from a slightly more general ``uncertainty
principle'' proven in \cite{Havin}. It concerns M. Rietz's
potentials of \textit{charges} (not necessarily absolutely
continuous with respect to Lebesgue measure) and $\alpha$'s not
necessarily from $(0,1)$; for the history of the problem and its connections with the uniqueness
problem for Laplase's equation see also \cite{HJ} and \cite{BH}. Havin \cite{Havin}
posted the following question: is it possible to omit H\"{o}lder's
condition on $f$ near $V$ in Theorem 1? Moreover, it was still
unknown if there exist a nonzero \textit{continuous} function
$f\in dom\, U_\alpha$ and a set $V$ of positive length satisfying
(\ref{allnull}).

In \cite{BH}, it was shown that the answer to the last question
is affirmative. However, the function $f$ constructed in
\cite{BH}, being continuous, does not satisfy any H\"{o}lder's
condition.

In this paper, we build a nonzero \textit{h\"{o}lderian} function
$f\in dom U_\alpha$ vanishing with its M. Rietz potential
$U_\alpha f$ on some set of positive length.

As well as in \cite{BH}, we build the desired function using
the techniques of ''correction`` which were proposed by D. Menshov and 
applied to the problems of potential theory in \cite{AK},
\cite{Binder}, \cite{BW}, \cite{wolff}. Our progress (as compared to \cite{BH}) 
is based on improving the correction
process using elementary probabilistic techniques (see Lemma
5\footnote{I am grateful to S. Smirnov for useful discussions
concerning this lemma.}). We also mention that we deal with
\textit{complex} (not only real-valued!) densities $f$. This
detail seems unessential, but it surprisingly simplifies the
construction of a desired (\textit{real-valued!}) density $f$, even
in case \cite{BH}, when the goal was just to find
\textit{continuous} (not necessarily  h\"{o}lderian) $f$.

The main result of this paper is the following theorem.
\begin{Theorem} There exist a nonzero function $f\in
dom\:U_\alpha$, a set $V\subset {\mathbb R}$ of positive length
and a positive number $r$, such that $f$ and $U_\alpha f$ are
identically zeroes on $V$ and $f$ satisfies H\"{o}lder's condition
with exponent $r$.
\end{Theorem}

{\bf Remark.}  The function $f$ we build is complex-valued.
In order to get a real-valued density, one should just take its
real or imaginary part; at least one of them will not be an
identical zero.

I am grateful to V. Havin for introducing me to the problem and
for useful discussions.

\bigskip \noindent{\bf1. Operator $W_\alpha$.} We will need an
operator that is in some sense inverse to $U_\alpha$. Let for
$g \in C^\infty_0({\mathbb R})$
$$
(V_\alpha g)(t):=\frac1\alpha\int_{{{\mathbb
R}}}{g(x)\frac{sgn(t-x)}{|t-x|^ \alpha}dx},
$$

$$
(W_\alpha g)(t)=((V_\alpha g)(t))',\quad t\in \mathbb{R}.
$$

Let $g$ be a function defined on $\mathbb{R}$, $\lambda>0$,
$\varepsilon>0$. Denote $(C_{\lambda}g)(x)=g(\lambda x)$,
$g_\varepsilon=\frac1\varepsilon C_{1/\varepsilon} g$. The next
lemma states some properties of $W$.
\begin{lemma}
\begin{enumerate}
\item $W_\alpha(C_0^\infty( {\mathbb R}))\subset {\rm dom}
U_\alpha$; \item $U_\alpha W_\alpha g = cg,g \in (C_0^\infty(
{\mathbb R}))$; \item $W_\alpha C_\lambda=\lambda^\alpha C_\lambda
W_\alpha$; \item $\alpha W_\alpha g = -g'\ast\frac{sgn
x}{|x|^\alpha}$ (we use $\ast$ for the convolution on $\mathbb{R}$);
\item $(W_\alpha g)(t)=(g\ast |x|^{-\beta)}(t), g \in {\rm dom}
W_\alpha, t \notin {\rm supp}\, g.$
\end{enumerate}
\end{lemma}
Hereinafter $\beta=\alpha+1$. All statements of this lemma
are well-known (2) or obvious (1,3,4,5). See, for example,
\cite[page 226]{BH}.

We use the following notation: $I=(-\frac12,\frac12)$; if $Q$ is a
bounded interval, then $c_Q$ is its center.

We write $\phi (t)$ for ``the finitizator'': $\phi  \in C^\infty$,
${\rm supp} \: \phi \subset I$, $\int_I \phi = 1$, $\phi \geq 0$.

In the proof we shall fix positive numbers $p$ and $\lambda$. For a function $h:\mathbb{R}\rightarrow\mathbb{C}$ 
we introduce its ``embedding to the interval $Q$''
$h_Q(t):=h(\frac{t-c_Q}{|Q|\lambda})$, $t\in\mathbb{R}$. Finally,
let $M_Q(h):=(\frac1{|Q|}\int_Q{|h|^p})^{1/p}.$

\bigskip\noindent{\bf2. Main Lemma.} The mail tool
for the proof of Theorem 1 will be Lemma 4. First, we prove
auxiliary lemmas 2 and 3. We state the existence of functions with
certain concrete numerical properties. These functions will serve
as ``building blocks'' for our construction. The  following proposition is principal for us: 
for some constant $B<0$ and some
$p>0$ there exists a function $h\in C_0^\infty$ with arbitrarily
small support satisfying $\int_{\mathbb R}(|1-W_\alpha h|^p-1)<B$.
The meaning of this fact is that it allows us to control both the
length of support of $h$ (future correcting term) and its influence on the ``amount'' of
potential $W_\alpha$ of the function we are correcting. We shall
rearrange this statement for it to be convenient for our purposes.

The function $h$ will be made from ``the finitizator'' $\phi$ by
means of an appropriate scaling.

Let for $\varepsilon>0$, $t\in {\mathbb R}$
$$F^{[\varepsilon]}(t):=(W_\alpha \phi_\varepsilon)(t), \quad
F^{[0]}(t):=|t|^{-\beta}.$$

\begin{lemma}
If a positive number $p$ is sufficiently small, then
$$
J(p) := \int_{\mathbb R}(|1-F^{[0]}|^p-1)dx < 0
$$
\end{lemma}

{\bf Proof.} Let $L:=\int_{\mathbb R}\log|1-F^{[0]}|$. Then, as
$(a^p-1)/p$ is monotone in $p$ for any $a>0$ and converges to
$\log a$ as $p\rightarrow 0$, we have $\lim_{p\searrow
0}\frac{J(p)}{p} = L$ (note that $||1-F^{[0]}(t)|^p-1| \leq
c|t|^{\beta}$, when $|t|$ is large, and if $p<1/\beta$, then the
integral $J(p)$ converges in zero as well). But 
$L$ can be computed exactly: $L=2\pi\cot\frac\pi\beta<0$. One can find the computation, 
for example, in \cite[page 233]{BH}.

From now on the number $p$ found in the previous lemma will be
fixed.

In the next lemma we pass from $F^{[0]}$ (the potenial $W_\alpha$
of the delta-function) to the potential of some concrete function
$\phi_\varepsilon$. We also introduce a ``small complex rotation'':
multiply $\phi_\varepsilon$ by $e^{i\theta}$ with small $\theta$.
This leads to some technical simplifications in the future.

\begin{lemma}
There exist numbers $B<0$, $\theta_0>0$ and $\varepsilon_0>0$,
such that if $0\leq\theta\leq\theta_0$ and
$0\leq\varepsilon<\varepsilon_0$, then
$$
J(\varepsilon,\theta):=\int_{\mathbb
R}(|1-e^{i\theta}F^{[\varepsilon]}|^p-1)<B.
$$
\end{lemma}
{\bf Proof.} Using the homogeneity property of $W_\alpha$ (point 3 of
Lemma 1), we get
\begin{equation}
\label{estimF} |F^{[\varepsilon]}(t)| \leq
C(\alpha)\min(\frac1{\varepsilon^\beta},\frac1{|t|^\beta}),\quad
t\in {\mathbb R}.
\end{equation}
It is clear that $F^{[\varepsilon]}$ converges to $F^{[0]}$
pointwise as $\varepsilon\rightarrow 0$. It follows from Lebesgue
dominated convergence theorem that
$$
\lim_{\varepsilon\searrow 0}J(\varepsilon,0)=J(p).
$$
Choose $B<0$ to satisfy $J(\varepsilon,0)<2B$ for
$\varepsilon\in(0,\varepsilon_0)$. We are to prove that
$J(\varepsilon,\theta)\stackrel{\theta\rightarrow
0}{\longrightarrow} J(\varepsilon,0)$ uniformly in $\varepsilon$.
Indeed,
$$
|J(\varepsilon,\theta)-J(\varepsilon,0)| = \left|\int_{\mathbb
R}|1-e^{i\theta}F^{[\varepsilon]}|^p-|1-F^{[\varepsilon]}|^p\right|
\leq
\left|\int_{|x|>\theta^{-p/2}}\right|+\left|\int_{|x|<\theta^{-p/2}}\right|=:J_1+J_2.
$$
We have
$$
J_1=\left|\int_{|x|>\theta^{-p/2}}|1-e^{i\theta}F^{[\varepsilon]}|^p-|1-F^{[\varepsilon]}|^p\right|
\leq 2C \left|\int_{|x|>\theta^{-p/2}}|t|^{-\beta}\right|
\stackrel{\theta\rightarrow 0}{\longrightarrow} 0,
$$
as the integral converges in infinity (and does not depend on
$\varepsilon$). Then, using the inequality $|a^p-b^p|\leq|a-b|^p$
for $a>0$, $b>0$ and $p\in(0,1)$, we get
$$
J_2=\left|\int_{|x|<\theta^{-p/2}}|1-e^{i\theta}F^{[\varepsilon]}|^p-|1-F^{[\varepsilon]}|^p\right|=
$$
$$
=\left|\int_{|x|<\theta^{-p/2}}|e^{-i\theta}-F^{[\varepsilon]}|^p-|1-F^{[\varepsilon]}|^p\right|
\leq
$$
$$
\leq \int_{|x|<\theta^{-p/2}}|e^{-i\theta}-1|^p \leq
\int_{|x|<\theta^{-p/2}}\theta^p\stackrel{\theta\rightarrow
0}{\longrightarrow} 0.
$$
So, $J_1+J_2\stackrel{\theta\rightarrow 0}{\longrightarrow} 0$
uniformly in $\varepsilon\in(0,\varepsilon_0)$, and the lemma is
proved.

We are now ready to formulate Lemma 4 -- the main tool of further
construction. We replace the whole line (from Lemma 3) by a
bounded interval, the constant 1 by an arbitrary function with a
small oscillation, and finally we bring in a (small) positive
parameter $\lambda$. This one is responsible for smallness of
the potential $W_\alpha$ of the correcting term far away from the
interval we correct on. Denote
$\gamma(\lambda):=(1+B\lambda/2)^{1/p}$, where $B$ is a constant
from Lemma 3. Notice that $0<\gamma(\lambda)<1$ for any sufficiently
small positive $\lambda$.

If $f$ is a function defined on some interval $Q$, we define
$osc_Q f:=\sup\limits_{x,y\in Q}(|f(x)-f(y)|)$ (\textit{the
oscillation} of $f$ on $Q$).

\begin{lemma}
There exist numbers $\theta>0$, $\lambda_0>0$, $\varepsilon_0>0$,
such that for any positive $\lambda< \lambda_0$ one can find a
number $\kappa>0$ with the following property: if $0<\varepsilon<
\varepsilon_0$, $Q$ is any bounded interval and a continuous
complex-valued function $h$ satisfies
\begin{equation}
\label{osc} osc_Q h \leq \kappa |h(c_Q)|,
\end{equation}
then

\begin{enumerate}
\item $M_Q(h-h(c_Q)e^{i\theta}F_Q^{[\varepsilon]}) \leq
\gamma(\lambda)|h(c_Q)|$; \item $\frac\theta2 |h(t)| \leq
|h(t)-h(c_Q)e^{i\theta}F_Q^{[\varepsilon]}(t)| \leq
\frac{C}{\varepsilon^\beta}|h(t)|,\quad t\in Q$.
\end{enumerate}
\end{lemma}

Recall that
$F^{[\varepsilon]}_Q(t)=F^{[\varepsilon]}(\frac{t-c_Q}{|Q|\lambda})$.

{\bf Proof.} Take $\varepsilon_0$ and $\theta$ from Lemma 3. First
we get an estimate in point 2:
$$
|h(t)-h(c_Q)e^{i\theta}F_Q^{[\varepsilon]}|\geq
|h(c_Q)||1-e^{i\theta}F_Q^{[\varepsilon]}(t)|-|h(t)-h(c_Q)|\geq
|h(c_Q)|\frac{3\theta}{4} - \kappa|h(c_Q)|=
$$
$$
=|h(c_Q)|(\frac{3\theta}{4} - \kappa)\geq
|h(t)|\frac1{1+\kappa}(\frac{3\theta}{4} - \kappa)\geq
|h(t)|\frac\theta{2}
$$
for any sufficiently small $\kappa$. We have used an elementary
inequality $dist(1,\{re^{i\theta}: r \in {\mathbb
R}\})=\sin\theta\geq\frac{3\theta}{4}$, if $\theta>0$ is small.

The right-hand inequality in the point 2 follows clearly from
(\ref{estimF}).

Now we prove point 1:
$$
(M_Q (h-h(c_Q)e^{i\theta}F_Q^{[\varepsilon]}))^p \leq (M_Q
(h-h(c_Q)))^p+|h(c_Q)|^p(M_Q(1-e^{i\theta}F_Q^{[\varepsilon]}))^p
\leq
$$
$$
\leq (osc_Q h)^p+|h(c_Q)|^p
(M_{\lambda^{-1}I}(1-e^{i\theta}F^{[\varepsilon]}))^p.
$$

Notice that
$$
(M_{\lambda^{-1}I}(1-e^{i\theta}F^{[\varepsilon]}))^p=
1+\lambda\int_{\lambda^{-1}I}(|1-e^{i\theta}F^{[\varepsilon]}|^p-1)=
$$
$$
=1+\lambda\int_{{\mathbb
R}}(|1-e^{i\theta}F^{[\varepsilon]}|^p-1)- \lambda\int_{{\mathbb
R} \backslash
\lambda^{-1}I}(|1-e^{i\theta}F^{[\varepsilon]}|^p-1).
$$
Taking (\ref{estimF}) into account, we get
$$
\left|\lambda\int_{{\mathbb R} \backslash
\lambda^{-1}I}(|1-e^{i\theta}F^{[\varepsilon]}|^p-1)\right| \leq
C\lambda^{\alpha+1}=o(\lambda),\lambda\longrightarrow 0,
$$
so, if $\lambda$ is sufficiently small, we have
$M_{\lambda^{-1}I}(1-e^{i\theta}F^{[\varepsilon]})\leq
1+B\lambda/2<1$. Therefore
$$
M_Q(h-h(c_Q)e^{i\theta}F^{[\varepsilon]})^p \leq
|h(c_Q)|^p(1+\left(\frac{osc_Q(h)}{h(c_Q)}\right)^p+2B\lambda/3)\leq
|h(c_Q)|^p(1+\kappa^p+2B\lambda/3)).
$$
If $\kappa$ is sufficiently small, then
$(1+\kappa^p+2B\lambda/3))<1+B\lambda/2$. The lemma is proved.

{\bf Remark 1.} Careful examination of the proof shows that one
can take $\kappa$ equal to
$\min((\frac{|B|\lambda}2)^{1/p},\frac\theta{8})$, if $\theta$ is
not too large.

{\bf Remark 2. } It is the left-hand inequality in the first point for
what we pass to complex-valued functions. One cannot obtain such
an estimate for real-valued functions.

{\bf Remark 3.} By this moment we have fixed parameters $p$ and
$\theta$. In what follows, constants that we regard as depending on $\alpha$
may also depend on these parameters. Later we shall fix an
appropriate $\lambda$ and, thus, $\kappa$ and $\gamma$.

\bigskip\noindent{\bf3. General idea of the construction.} Now we describe the 
plan of construction of $f$ and $V$ (see statement of Theorem 1). We shall build a sequence of functions $g_n$,
$g_n=g_{n-1}-r_{n-1}$, and a decreasing sequence of sets $V_n\subset
I$ with the following properties:
\begin{enumerate}
\item A nonzero function $g_1$ belongs to $C_0^\infty$, and
$supp\;g_1\subset\mathbb{R}\backslash I$; \item $r_k\in
C_0^\infty$, and $supp\;r_k\subset I$ for all $k\in\mathbb{N}$;
\item $\sum\limits_{k=1}^\infty|supp\;r_k|<\frac14$; \item
$|\bigcap\limits_{k=1}^\infty V_n|>\frac34$; \item
$\int_{V_n}|f_n|^p\stackrel{n\rightarrow\infty}{\longrightarrow}0$;
here $f_n:=W_\alpha g_n$ and $p$ is the positive number fixed in
the previous section \item Sequences $g_n$ and $f_n$ converge
uniformly on $\mathbb{R}$ to some continuous functions $g$ and $f$ correspondingly, and
$g=U_\alpha f$.
\end{enumerate}

Let $V:=\bigcap\limits_{k=1}^\infty V_n,$ $V':=\{x\in
I:\;g(x)=0\}.$ It follows from 1, 3 and 4 that $|V'|>\frac34$ and
$|V|>\frac34.$ Therefore $|V\cap V'|>\frac12$. It follows from 5
and 6 that $f|_V=0$. Finally, using 2, we conclude that
$g|_{\mathbb{R}\backslash I}=g_1|_{\mathbb{R}\backslash I}$, and,
thus, the function $g$ is not identically 0. Hence the set $V\cap
V'$ and the function $f$ satisfy all conditions of Theorem 1,
except (may be) H\"{o}lder's condition.

Now we describe more precisely the structure of sets $V_n$ and
correcting terms $r_k$. Bring in a {\sl sequence of positive
numbers $\{\delta_n\}_1^\infty$}. It will have the following
properties: $\delta_1=1, \frac{\delta_n}{\delta_{n+1}} \in
{\mathbb N}$. We denote the partition of the interval $I$ to
intervals of length $\delta_n$ as $H_n$. The set $V_n$ will be
obtained as the union $\bigcup_{Q\in G_n} Q$, where $G_n$ is some
subset of $H_n$. Roughly speaking, the set $G_n$ consists of all
intervals on which we haven't finish correction yet, in
particular, for all $k>n$ there holds $supp\;r_k\subset V_n$.

Let us fix a sequence of positive numbers
$\{\varepsilon_n\}_{n=1}^\infty$, such that
$\sum\limits_{n=1}^\infty\varepsilon_n<\frac14$, and, in addition,
$\varepsilon_n$ decay not very fast: $\varepsilon_n^{-1}=O(n^m)$
for some $m>0$. It will be responsible for the length of supports
of $r_n$: there will hold an estimate
$|supp\;r_n|\leq\varepsilon_n$ for all $n$.

We also demand $supp\;g_1\in(\frac12,\frac32)$, and, moreover,
$\forall t\in I f_1(t)\neq 0$. One can take, for
example, $g_1:=\phi(x-1)$.

Then, we choose some subset $G_n^g\subset G_n$ and let
\begin{equation}
\label{defr_n} r_n:=\sum\limits_{Q\in
G^g_{n+1}}(\lambda\delta_{n+1})^\alpha
f_n(c_Q)(\phi_{\varepsilon_n})_Q e^{i\theta}.
\end{equation}
So,
$$
W_\alpha r_n=\sum\limits_{Q\in G^g_{n+1}}
f_n(c_Q)F^{[\varepsilon_n]}_Q e^{i\theta}.
$$
Note that in such a definition condition 3 will be
satisfied by choosing the sequence $\varepsilon_n$ as described abowe.

The idea is that if $\delta_{n+1}$ is sufficiently small, then on
each interval $Q\in G^g_{n+1}$ the oscillation of $f_n$ is small
(there holds estimate (\ref{osc})), and one can apply Lemma 4
with $f_n$ as $h$. Its result, together with an observation that
functions $F^{[\varepsilon]}_Q$ decay sufficiently fast far away from
$Q$, allow us, using the notation $V_n^g:=\bigcup\limits_{Q\in G_n^g} Q$, to
prove an estimate
\begin{equation} \label{lessEta}
\int_{V^g_{n+1}}|f_{n+1}|^p \leq \eta \int_{V^g_n}|f_n|^p
\end{equation}
with some $\eta\in(0,1)$. If one chooses $G^g_n$ appropriately (if
on each step they occupy a large part of $G_n$), this leads to
an estimate of integral over the whole set $V_n$:
\begin{equation}
\label{lesseta2} \int_{V_{n}}|f_{n}|^p = O(\eta^{\frac{n}{2}}).
\end{equation}
Hence condition 5 will be obtained.

{\bf Remark 1.} The choice of $G_n$ (decreasing $V_n$ on each
step) allows us to make functions $f_n$ converge not only in the
sense of $L^p(I)$, but uniformly, in partiular, we get an estimate
$|f_n(c_Q)|=O(\eta'^n),\;Q\in G_n$, where $\eta'\in(0,1)$.

{\bf Remark 2.} If we did not worry about the control over
modulus of continuity, we could take $G_n^g:=G_n$. Then
(\ref{lesseta2}) automatically follows from (\ref{lessEta}), and
the whole construction becomes more simple. Unfortunately, in order to
get H\"{o}lder's condition, one should pick out the set $G_n^g$ on
each step -- this is a set of intervals where the oscillation of $f_n$
is especially small -- and make the correction only there.

\bigskip\noindent{\bf4. Remarks on the estimate of modulus of continuity.}
In this section, we explain (not quite rigorously) what does
estimates of modulus of continuity of $f$ depend on.

We use the following simple fact: \emph{if a sequence of functions
$h_n$ converges on $\mathbb{R}$ to the function $h$, whereas
$|h_n-h|\leq C_1\eta_1^n$, and $|h_n'|\leq C_2 R^n$ (here
$\eta_1\in(0,1),\; R>1$), then $h$ satisfies H\"{o}lder's
condition with an exponent $\log\eta_1/\log\frac{\eta_1}{R} $}.

The condition $|f_n-f|\leq C_1\eta_1^n$, will follow from Remark 1 at 
the end of the previous section (and, in fact, from estimate (\ref{lesseta2})). 
When one estimates the derivative
$f_n'(t)$ of the function $f_n$, the main role is played by the
last added term $W_\alpha r_{n-1}$, or, more precisely,
the building block $f_{n-1}(c_{Q_t})F^{[\varepsilon_{n-1}]}_{Q_t}
e^{i\theta}$, where the interval $Q_t\in H_n$ is defined by the
statement $t\in Q_t$. From the homogeneity properties of $W_\alpha$ one
can get (for $Q\in H_n$) an estimate

\begin{equation} \label{estimFP}
|(F^{[\varepsilon_{n-1}]}_{Q})'(t)|\leq
\frac{c}{\varepsilon^{-\beta-1}\lambda\delta_n}\;t\in\mathbb{R}
\end{equation}

Thus, one can get (say, on $V_n$) an estimate

\begin{equation} \label{ocPrime1}|f_n'|\leq
\frac{Cf_n(c_{Q_t})}{\delta_n\varepsilon_{n-1}^{\beta+1}},
\end{equation}

and so, everywhere,

\begin{equation} \label{ocPrime}|f_n'|\leq
\frac{C\eta'^n}{\delta_n\varepsilon_{n-1}^{\beta+1}}.
\end{equation}

This means that the obtained function $f$ would be h\"{o}lderian
if numbers $\delta_n^{-1}$ did not grow faster than some geometric
series, in other words, if every time we divided the interval $Q\in
H_n$ into the same number of parts. On the other hand, the exponent
$\log\eta_1/\log\frac{\eta_1}{R}$ tends to zero as $R\rightarrow
\infty$, hence it is clear that if
$\frac{\delta_n}{\delta_{n+1}}\rightarrow \infty$ as
$n\rightarrow\infty$, then we are unable to prove H\"{o}lder's
condition with any exponent.

It is clear, however, that if in order to define $\delta_{n+1}$ we
use a natural estimate (\ref{ocPrime1}) (recall what the necessity
to choose small $\delta_{n+1}$ is due to: we need to estimate the
oscillation of $f_n$ in order to use Lemma 4 -- condition
(\ref{osc})), then because of the increasing multiplier
$\varepsilon_{n-1}^{-\beta-1}$ we should take
$\delta_{n+1}/\delta_n$ tending to zero to make (\ref{osc}) hold.
Therefore we need finer esimates of modulus of continuity of
$f_n$, holding, however, not on the whole set $V_{n+1}$, but on
some ''good`` part $V_{n+1}^g$ of it.

Note that the building block $F^{[\varepsilon_n]}_Q$ and its
derivative $(F^{[\varepsilon_n]}_Q)'$ are large in modulus (as
$\varepsilon_n^{-\beta}$ and
$\delta_{n+1}^{-1}\varepsilon_n^{-\beta-1}$ correspondingly) only
near the center of $Q$; if we consider them outside the interval
of the length $\tau |Q|$ and with the same center with $Q$, where
$0<\tau<1$, then we get estimates $|F^{[\varepsilon_n]}_Q|\leq C$
and $|(F^{[\varepsilon_n]}_Q)'|\leq C\delta_{n+1}^{-1}$.

One has an idea - exclude this ``bad'' central part of $Q$ and
correct on the remaining part only. Unfortunately, if we drop it
forever, this will mean that on each step one eliminates from
$V_n$ a subset of length $\tau|V_n|$, and this makes $\cap_{n\in
{\mathbb N}}V_n$ have zero length.

Therefore, for each interval from $H_n$ we bring in a system of
its ``bad'' subsets, and on each of them we shall ``make a pause''
- not correct during the next few steps, until the partition
$H_{n+k}$ becomes so fine that the estimates of $f_n$ and its
derivative become satisfactory (requirements to this estimates
becomes weaker if $n$ grows). The the pause duration depends on
the distance of the corresponding subset from the center of the
interval, i.e. on how ``bad'' $f_n$ is on this subset. According
to this, $G_{n+1}$ is divided into two parts: $G_{n+1}^g$, where
the correction is made on this step, and $G_{n+1}^d$, where we
do not do anything for the time being.\footnote{Top indices $g$
and $d$ are the first letters of ``go'' and ``delay''}

After that, we shall estimate an amount of intervals from $G_n$
such that have more than one half of their ``ancestors'' from $G_k,\;
k=0,1,\dots n-1,$, belong to $G_k^d$. It turns out that there
are only a few of them (if $G_k^d$ is a small part of $G_k$ for
each $k$), and we drop them out. For the rest, we prove an
estimate like (\ref{lesseta2}), using Lemma 4 and fast decay of
$F^{[\varepsilon_{n}]}_{Q}$ away from $Q$.

\bigskip\noindent{\bf5. Definition of sets $G_n^g$.} To complete the
construction, we should determine the sequence $\delta_n$ and sets
$G_n$ and $G_n^g$. We need a number of estimates depending on how
one picks out ``good'' subsets $G_n^g$ from $G_n$, but not on
the way to choose $G_n$ themselves. We prove these estimates
in sections 6, 7 and 8. Later, in section 9, we will define the way
to choose $G_n$.

Bring in a positive parameter $\delta$, such that $\delta^{-1}\in
{\mathbb N}$, and let $\delta_n:=\delta^n$.

Then, bring in a parameter $\tau\in (0,1)$. Let
$\tau^{-1}\in{\mathbb N}$, and, moreover, $\tau/\delta\in{\mathbb
N}$. We use the following notation: if $a>0$ and $Q$ is a bounded
interval, then $Q[a]:=Q\backslash Q'$, where $Q'$ is an interval
of the length $a|Q|$ and with the same center as $Q$.

We are now ready to define the set $G_{n+1}^g$. An interval $Q\in
G_{n+1}$ belongs to $G_{n+1}^g$, iff for any $k=0,1,\dots, n-1$
there holds an implication $Q\subset Q'\in G_{n-k}^g\Rightarrow
Q\subset Q'[\tau^{k+1}]$. In other words, if in the $(n-k)$-th step we
have made a correction\footnote{Recall that it means that $Q'$
belongs to the set of indices of summation in the definition
(\ref{defr_n}) of the corresponding correcting term $r_{n-k-1}$} on
the interval $Q'$, then on the next step the correction is
forbidden on the set $Q'\backslash Q'[\tau]$, on the $(n-k+2)$-th step it
is forbidden on the set $Q'\backslash Q'[\tau^2]$, and so on. It
follows from the condition $\tau/\delta\in{\mathbb N}$ that the
interval $Q\in G_{n+1},\;Q\subset Q'\in G_{n-k}^g$ either lies in
$Q'[\tau^{k+1}]$ or does not intersect it.

In fact, $Q'[\tau^{k+1}]\backslash Q'[\tau^k]$ are the very
``bad'' subsets of $Q'$; on the $k$-th of them the ``length of pause''
is $k$ steps.

Let us make some simple, but important observations. It is
easy to see that if $Q$ is a bounded interval and
$dist(t,c_Q)>3\lambda\varepsilon_n|Q|$, then for any $n\in
{\mathbb N}$
\begin{equation}
\label{Fdaleko} |F^{[\varepsilon_n]}_Q(t)|\leq
C(\alpha)\frac{(\lambda|Q|)^\beta}{|t-c_Q|^\beta},
\end{equation}
and, besides,
\begin{equation}
\label{FPdaleko} |(F^{[\varepsilon_n]}_Q)'(t)|\leq
C(\alpha)\frac{(\lambda|Q|)^\beta}{|t-c_Q|^{\beta+1}}.
\end{equation}
This means that for all $k\in {\mathbb N}$ and for $t\in
Q[\tau^k]$
\begin{equation}
|F^{[\varepsilon_n]}_Q(t)|\leq
C_1(\alpha)\frac{\lambda^\beta}{\tau^{k\beta}}, \label{goodF}
\end{equation}
\begin{equation}
|(F^{[\varepsilon_n]}_Q)'(t)|\leq
C_1(\alpha)\frac{\lambda^\beta}{|Q|\tau^{k(\beta+1)}}.
\label{goodPrime}
\end{equation}
Indeed, for $\tau^k/2>3\lambda\varepsilon_n$ these estimates
coincide with the previous ones, and for $\tau^k/2>3\lambda\varepsilon_n$
we can use estimates
$$
|F^{[\varepsilon_n]}_Q(t)|\leq
\frac{C(\alpha)}{\varepsilon_n^\beta}
$$
and
$$
|(F^{[\varepsilon_n]}_Q)'(t)|\leq
\frac{C(\alpha)}{|Q|\lambda\varepsilon_n^\beta},\quad
t\in\mathbb{R}.
$$
We can improve the right-hand inequality in point 2 of Lemma 4 for
$t\in Q[\tau^k]$. Indeed, applying (\ref{goodF}) instead of
(\ref{estimF}), we get
\begin{equation}
|h(t)-h(c_Q)e^{i\theta}F_Q^{[\varepsilon]}(t)| \leq
\frac{C(\alpha)\lambda^\beta}{\tau^{k\beta}}|h(t)|,\quad t\in
Q[\tau^k].
\end{equation}
Let $t\in I$. Denote as $Q^n_t$ an element of $H_n$ defined by the
condition $t\in Q^n_t$. Let $D^k_n(t):=\#\{l=1,\dots,n: t\in
Q^l_t\backslash Q^l_t[\tau^k]\}$. For $t\in V_n$ let
$\widetilde{D}_n(t):=\#\{l=1,\dots,n: Q^l_t\in G^d_l\}$. We know
that if $Q^l_t\in G^d_l$, then for some $k_l\in {\mathbb N}$ there
holds an inclusion $Q^l_t\subset Q^{l-k_l}_t\backslash
Q^{l-k_l}_t[\tau^{k_l}]$. Moreover, for $l_1\neq l_2$ either
$k_{l_1}\neq k_{l_2}$, or $l_1-k_{l_1}\neq l_2-k_{l_2}$. So, with
each natural $l\leq n$, such that $Q^l_t\in G^d_l$, we can
associate a pair $(l-k_l,k_l)$, such that $Q^l_t\subset
Q^{l-k_l}_t\backslash Q^{l-k_l}_t[\tau^{k_l}]$, and such a mapping
will be injective. Hence $\widetilde{D}_n(t)\leq \sum\limits_{k\in
{\mathbb N}}D_n^k(t)=:D_n(t)$. The next section is devoted to
estimates of lengths of sets $E_n:=\{t\in I:\: D_n(t)\geq n/2\}$.

\bigskip\noindent{\bf6. Estimates of lengths of sets $E_n$.} \begin{lemma}. For $\tau$
sufficiently small there holds an inequality
\begin{equation}
|E_n|\leq \frac{C(\tau)}{n^2},
\end{equation}
where $C(\tau)$ tends to zero as $\tau\rightarrow 0$.
\end{lemma}
{\bf Proof.} Let
$$
{\xi'}_i^{(k)}:=\sum\limits_{Q\in H_i}\chi_{Q\backslash
Q[\tau^k]};
$$
$$
\xi_i^{(k)}:={\xi'}_i^{(k)}-\tau^k.
$$
Then, considering $\xi_i^{(k)}$ as random variables on the
probabilistic space $I$ with the measure $dx$, we have
$E\xi_i^{(k)}=0$. We can describe sets $E_n$ in terms of functions
$\xi_i^{(k)}$ as follows:
$$
x \in E_n \Leftrightarrow \sum_{k=1}^{n-1} \sum_{i=1}^{n-k}
{\xi'}_i^{(k)}(x) \geq \frac{n}{2}
$$
Hence we are to estimate probabilities of the event that sums
of random variables $\sum_{i=1}^{n-k} {\xi'}_i^{(k)}$ are large.
It is easy to see that random variables ${\xi'}_i^{(1)},\quad
i=1,2,\dots$ are independent: it follows from the fact that
$\tau/\delta$ is an integer. Unfortunately, for $k>1$ one cannot say
the same thing about ${\xi'}_i^{(k)},\quad i=1,2,\dots$, because
$\tau^k/\delta$ is not necessarily an integer. But, for $k>1$
variables ${\xi'}_i^{(k)},\quad i=1,2,\dots$ are still in some
sense ``almost independent'', and we use it.

Note that if $j\geq i+k$, then $\xi_i^{(k)}$ is a constant on each
interval from $H_j$ (it follows from inclusion
$\tau^k/\delta^k\in\mathbb{N}$). Hence we can made the following
observation:
\begin{enumerate}\item If $i_1+k\leq i_2\leq i_3\leq i_4$, then
$E(\xi^{(k)}_{i_1}\xi^{(k)}_{i_2}\xi^{(k)}_{i_3}\xi^{(k)}_{i_4})=0$
(as the function $\xi^{(k)}_{i_2}\xi^{(k)}_{i_3}\xi^{(k)}_{i_4}$
is periodic with a period equal to $\delta^{i_1+k}$, and
$\xi^{(k)}_{i_1}$ is constant on each interval from $H_{i_1+k}$).
\item If $i_1\leq i_2\leq i_3\leq i_4-k$, then
$E(\xi^{(k)}_{i_1}\xi^{(k)}_{i_2}\xi^{(k)}_{i_3}\xi^{(k)}_{i_4})=0$
(as the function $\xi^{(k)}_{i_4}$ is periodic with a period equal
to $\delta^{i_3+k}$, and
$\xi^{(k)}_{i_1}\xi^{(k)}_{i_2}\xi^{(k)}_{i_3}$ is constant on
each interval from $H_{i_3+k}$).
\end{enumerate}
Now we write
$$ P(|\sum_{i=1}^n\xi_i^{(k)}|> \varepsilon)\leq
\frac{E(\xi_{1}^{(k)}+\dots+\xi_{n}^{(k)})^4}{\varepsilon^4}=$$
\begin{equation}=
\frac{\sum\limits_{(i_1,i_2,i_3,i_4)\in\{1,\dots,n\}^4}
E(\xi_{i_1}^{(k)}\xi_{i_2}^{(k)}\xi_{i_3}^{(k)}\xi_{i_4}^{(k)})}{\varepsilon^4}.
\label{sum}
\end{equation}
First note that
$E(\xi_{i_1}^{(k)}\xi_{i_2}^{(k)}\xi_{i_3}^{(k)}\xi_{i_4}^{(k)})\leq
E|\xi_{i_1}^{(k)}|=2\tau^k(1-\tau^k)\leq 2\tau^k$ (the first
inequality follows from the inequality $|\xi_{i}^{(k)}|<1$ for all
$i$ and $k$). Second, if $(j_1,j_2,j_3,j_4)$ is a non-decreasing
permutation of the numbers $(i_1,i_2,i_3,i_4)$, then a term
$E(\xi_{i_1}^{(k)}\xi_{i_2}^{(k)}\xi_{i_3}^{(k)}\xi_{i_4}^{(k)})$
may differ from zero only if $j_2-j_1<k$ and $j_4-j_3<k$ (by the
above observation). But the number of such fours
$(j_1,j_2,j_3,j_4)$ does not exceed $\frac{n(n-1)}2k^2$. Hence, the
number of nonzero terms in the numerator in (\ref{sum}) does not
exceed $4!\frac{n(n-1)}2k^2$. Therefore

\begin{equation}
P(|\sum_{i=1}^n\xi_i^{(k)}|>\varepsilon)\leq
\frac{Cn^2k^2\tau^k}{\varepsilon^4}. \label{cheb}
\end{equation}

Let $E_n^k:=\{t\in I: D^k_n(t)\geq
n(\sqrt[8]{(4\tau)^k}+\tau^k)\}$. If $\tau$ is so small that
$\sum_1^\infty\sqrt[8]{(4\tau)^k}+\tau^k<\frac12$, then
$|E_n|\leq\sum_{k=0}^\infty|E_n^k|$. Finally, note that
$D_n^k(t)=\sum_{i=1}^n{\xi'}^{(k)}_i(t)=\sum_{i=1}^n{\xi}^{(k)}_i(t)+n\tau^k$,
so, using (\ref{cheb}), we get
$$
|E_n|\leq\sum\limits_{k=1}^\infty|E_n^k|\leq
\sum\limits_{k=1}^\infty\frac{Ck^2\tau^\frac{k}{2}}{n^2 2^k}.
$$
The lemma follows from this estimate.

\bigskip\noindent{\bf7. Estimates of $f_n$.} In this section, we
prove some estimates for $f_n$'s, main of them are estimate
(\ref{kappa}), which allows us to apply Lemma 4, and estimate
(\ref{ghvost}), which shows that for $\lambda$ sufficiently small
the terms $f_n(c_{Q'}) F^{[\varepsilon_n]}_{Q'} e^{i\theta}$,
corresponding to $Q'\neq Q$, do not change the
situation on $Q$ essentially. All the estimates we prove do not depend on the
choice of $G_n$, but depend on how we pick out subsets $G_n^g$
(namely, we use estimates (\ref{goodF}) and (\ref{goodPrime})).
The sets $G_n$, as mentioned above, will be defined later. Let
\begin{equation}\label{defT}
T_{n+1}(t):=\sum\limits_{Q\in G^g_{n+1},Q\neq Q^{n+1}_t}f_n(c_Q)
F^{[\varepsilon_n]}_Q e^{i\theta}.
\end{equation}
Recall that an interval $Q_t^{n+1}\in H_{n+1}$ is defined by
condition $t\in Q_t^{n+1}$.

\begin{lemma} There exists a positive number $\rho=\rho(\alpha)$,
such that for $\lambda>0$ sufficiently small and for
$\delta=\delta(\lambda)>0$ sufficiently small there holds the
following:
\begin{enumerate}\item for all $t\in I$\begin{equation}|T_{n+1}(t)|\leq
\frac{c(\alpha)\lambda^\beta|f_n(t)|}{\rho};\label{ghvost}
\end{equation}\item if $k\leq n$, $x\in V_{n+1}^g$, $y\in I$ and
$|x-y|\leq \delta^k$, then
\begin{equation}
\label{gfxfy} |f_n(x)|\leq\frac{|f_n(y)|}{\rho^{n-k+1}};
\end{equation}
\item for all $t\in I$, there holds an estimate
\begin{equation}
|T'_{n+1}(t)|\leq
\frac{c(\alpha)\lambda^\beta|f_n(t)|}{\delta^{n+1}\rho};
\label{ghprime}
\end{equation}
\item for all $t\in V_{n+1}^g$, there holds an estimate
\begin{equation}
|f_n'(t)|\leq \frac{c_1(\alpha)\lambda^\beta|f_n(t)|}{\delta^n};
\label{fprime}
\end{equation}
\item for all $Q\in G_{n+1}^g$, there holds an estimate
\begin{equation}
osc_Qf_n\leq \kappa|f_n(c_Q)|\label{kappa}.
\end{equation}
\end{enumerate}
\end{lemma}

{\bf Proof.} Inequality (\ref{kappa}) clearly follows from
(\ref{fprime}), if we take $\delta$ sufficiently small (an
interval $Q$ in (\ref{kappa}) is of the length $\delta^{n+1}$). We
derive (\ref{ghvost}) and (\ref{ghprime}) from (\ref{gfxfy}),
which, in turn, follows from (\ref{ghvost}) and (\ref{kappa}) for
preceding n. Finally, (\ref{fprime}) follows from (\ref{kappa})
and (\ref{ghprime}) for preceding n.

The base of induction -- (\ref{kappa}) and (\ref{gfxfy}) for $n=1$
-- is provided by the condition $f_1(t)\neq 0,\;t\in I$ (see
Section 3) and the choice of $\rho$ and $\delta$ sufficiently
small.

Derive (\ref{ghvost}) and (\ref{ghprime}) from (\ref{gfxfy}). Fix
$t\in I$. Denote as $G_\epsilon$ the set of all intervals $Q'\in
H_{n+1}$ satisfying the property $dist\,(c_{Q'},Q_t^{n+1})\geq
\epsilon$. Denote
$$
\sigma_\epsilon(t):=\sum\limits_{Q'\in
G_\epsilon}|F^{[\varepsilon_n]}_{Q'}(t)|,
$$
$$ \sigma^*_\epsilon(t):=\sum\limits_{Q'\in
G_\epsilon}|(F^{[\varepsilon_n]}_{Q'})'(t)|. $$

We need estimates
\begin{equation} \label{sigma} \sigma_\epsilon(t) \leq
c(\alpha)\lambda^\beta\left(
\frac{\delta^{n+1}}\epsilon\right)^\alpha
\end{equation}
and
\begin{equation}\label{sigmaP} \sigma^*_\epsilon(t)\leq
c(\alpha)\lambda^\beta\frac{\delta^{(n+1)\beta}}{\epsilon^{\beta+1}}.
\end{equation}
One can obtain them by estimating each term by 
(\ref{Fdaleko}) and (\ref{FPdaleko}) correspondingly, and then estimating the
sum by an integral. The detailed proof of the first one can be
found in \cite[page 234]{BH}, the second one can be proved in the
same way.

In order to get the above estimate of $|T|$ and $|T'|$, we shall
divide terms in the right-hand side of (\ref{defT}) into several
groups according to their distance from the point $t$. For each
group, we estimate $|f(c_Q)|$ by means of (\ref{gfxfy}) (the closer to $t$ the interval $Q$ is, 
the better is this estimate) and then apply
(\ref{sigma}) (correspondingly, (\ref{sigmaP}) for $|T'|$), which,
by contrast, becomes better when $\epsilon$ grows.

So, let $G_{n+1}^g:=\bigsqcup_{k\leq n+1}G^{[k]}$, where
$G^{[n+1]}:=G^g_{n+1}\backslash G_{\delta^n}$,
$G^{[k]}:=(G^g_{n+1}\cup G_{\delta^k})\backslash
G_{\delta^{k-1}}$. For  $y\in G^{[k]}$ we have $|f_n(y)|\leq
|f_n(t)|/\rho^{n-k+2}$ and $dist(G^{[k]},t)\geq \delta^k/2$,
therefore
$$
|T_{n+1}(t)|\leq \sum_k\sum\limits_{Q'\in
G^{[k]}}|F^{[\varepsilon_n]}_{Q'}(t)||f_n(t)|/\rho^{n-k+2}\leq
$$
$$
\leq
c(\alpha)\lambda^\beta|f_n(t)|\sum_k\frac{\delta^{(n+1)\beta}}{\delta^{k\beta}\rho^{n-k+2}}
\leq
c(\alpha)\lambda^\beta|f_n(t)|\rho^{-1}\sum_k\left(\frac{\delta^\beta}{\rho}\right)^{n-k+1};
$$
similarly
$$
|T'_{n+1}(t)|\leq \sum_k\sum\limits_{Q'\in
G^{[k]}}|(F^{[\varepsilon_n]}_{Q'})'(t)||f_n(t)|/\rho^{n-k+2}\leq
$$
$$
\leq c(\alpha)\frac{\lambda^\beta|f_n(t)|}{\delta^{n+1}}
\sum_k\frac{\delta^{(n+1)(\beta+1)}}{\delta^{k{\beta+1}}\rho^{n-k+2}}
\leq c(\alpha)\frac{\lambda^\beta|f_n(t)|}{\delta^{n+1}\rho}\sum_k
\left(\frac{\delta^{\beta+1}}{\rho}\right)^{n-k+1}.
$$
Taking $\delta$ according to the condition
$\frac{\delta^\beta}\rho<\frac12$, we get (\ref{ghvost}) and
(\ref{ghprime}).

Now let us prove that (\ref{gfxfy}) and (\ref{fprime}) follow from (\ref{ghvost}), (\ref{ghprime}) and
(\ref{kappa}) for the previous
$n$. We need an estimate
\begin{equation} |f_n(x)|\leq
\frac{4}{\theta}|f_{n+1}(x)|\leq\dots\leq
\left(\frac4\theta\right)^k|f_{n+k}(x)|,\quad x\in I,\;n,k\in
{\mathbb N}, \label{snizu}
\end{equation}
which holds for $\lambda$ sufficiently small. Let us prove it. Let
$Q_x^{n+1}\in G_{n+1}^g$. Then we have
$$
|f_n(x)|\leq\frac2\theta|f_n(x)-f_n(c_{Q_x})e^{i\theta}F^{[\varepsilon_n]}_{Q_x}(x)|\leq
$$
$$
\leq
\frac2\theta(|f_{n+1}(x)|+|T_{n+1}(x)|)\leq\frac2\theta(|f_{n+1}(x)|+c(\alpha)\lambda^\beta\rho^{-1}|f_n(x)|).
$$
The first inequality follows from point 2 of Lemma 4 (which is
applicable because of (\ref{kappa})), the last one -- from
(\ref{ghvost}). Now, the inequality $|f_n(x)|\leq
\frac{4}{\theta}|f_{n+1}(x)|$ follows from the last estimate, if
$2c(\alpha)\lambda^\beta\theta^{-1}\rho^{-1}<\frac12$. If
$Q_x^{n+1}\notin G_{n+1}^g$, then it follows from (\ref{ghvost})
even easier. So, (\ref{snizu}) is proved.

Now, let $k\leq n+1$, $x\in V_{n+2}^g$, $y\in I$ and $|x-y|\leq
\delta^k$. Let $k':=\max\{l\leq n: x\in V_{l+1}^g\}$. From the
fact that $x$ is again in a ``good'' set $V_{n+2}^g$, it follows
that
\begin{equation}x\in
Q^{k'+1}_x[\tau^{n-k'+1}] \label{include}.
\end{equation}
First assume $k\leq k'$. Then
$$
|f_{n+1}(x)|=|f_n{x}|+|T_{n+1}(x)|\leq
|f_n(x)|(1+c(\alpha)\frac{\lambda^\beta}\rho)\leq\dots\leq
$$
$$
\leq |f_{k'+1}(x)|(1+c(\alpha)\frac{\lambda^\beta}\rho)^{n-k'}\leq
$$
$$
\leq (|f_{k'}(x)|+|T_{k'+1}(x)|+|f_{k'}(c_{Q_x^{k'}})|
|F^{[\varepsilon_n]}_{Q_x^{k'+1}}|)(1+c(\alpha)\frac{\lambda^\beta}\rho)^{n-k'}\leq
$$
$$
\leq |f_{k'}(x)|(1+c(\alpha)\frac{\lambda^\beta}\rho+(1+\kappa)
|F^{[\varepsilon_n]}_{Q_x^{k'+1}}|)(1+c(\alpha)\frac{\lambda^\beta}\rho)^{n-k'}.
$$
Applying (\ref{goodF}) and taking into account (\ref{include}), we
write on:
$$
|f_{n+1}(x)|\leq
|f_{k'}(x)|(1+C(\alpha)\lambda^\beta(\frac{1}\rho+\frac
{(1+\kappa)}{\tau^{(n-k'+1)\beta}})
)(1+c(\alpha)\frac{\lambda^\beta}\rho)^{n-k'}\leq
$$
$$
\leq
\frac{|f_{k'}(y)|}{\rho^{k'-k+1}}(1+C(\alpha)\lambda^\beta(\frac{1}\rho+\frac
{(1+\kappa)}{\tau^{(n-k'+1)\beta}})
)(1+c(\alpha)\frac{\lambda^\beta}\rho)^{n-k'}\leq
$$
$$
\leq |f_{n+1}(y)|\frac{\left(
\frac4\theta\right)^{n-k'+1}}{\rho^{k'-k+1}}(1+C(\alpha)\lambda^\beta(\frac{1}\rho+\frac
{(1+\kappa)}{\tau^{(n-k'+1)\beta}}))(1+c(\alpha)\frac{\lambda^\beta}\rho)^{n-k'}.
$$
We have used the induction assumption (\ref{gfxfy}) for $n=k'$ and
(\ref{snizu}). Now, if we choose $\rho$ so that
$\sqrt{\rho}<\theta/50$ and
$\sqrt{\rho}<\tau^\beta/50(C(\alpha))$, we get
$$|f_{n+1}(x)|\leq
|f_{n+1}(y)|\frac1{\rho^{k'-k+1}}5^{-n-k'+1}\rho^{-n+k'-1}\cdot
$$
$$
\cdot
((1+c(\alpha)\frac{\lambda^\beta}\rho)^{n-k'+1}+5^{-n-k'+1}\rho^{-n+k'-1}(1+c(\alpha)\frac{\lambda^\beta}\rho)^{n-k'}).
$$
Then we choose $\lambda=\lambda(\rho)$ so that
$(1+c(\alpha)\frac{\lambda^\beta}\rho)<2$, and get
$$
|f_{n+1}(x)|\leq \frac{|f_{n+1}(y)|}{\rho^{n-k'+1+k'-k+1}},
$$
and we are done.

Now we prove (\ref{fprime}). Let $t\in V^g_{n+2}$. Let, as in the
proof of (\ref{gfxfy}), $k':=\max\{l\leq n: x\in V_{l+1}^g\}$.
Again there holds (\ref{include}). We have:
$$
|(f_{n+1})'(t)|\leq|f'_n(t)|+|T'_{n+1}(t)|\leq\dots\leq
$$\begin{equation}
\leq|f'_{k'}(t)|+
|f_{k'}(c_{Q_t^{k'+1}})||(F^{[\varepsilon_{k'}]}_{Q_t^{k'+1}})'(t)|+
\sum\limits_{l=k'+1}^{n+1}|T'_{l}(t)|. \label{prproof}
\end{equation}
Applying the induction assumption ((\ref{fprime}) for $n=k'$) and
(\ref{snizu}), we have
$$
|f'_{k'}(t)|\leq
c_1(\alpha)\lambda^\beta\frac{|f_{k'}(t)|}{\delta^{k'}}\leq
c_1(\alpha)\lambda^\beta(4/\theta)^{n-k'+1}\frac{|f_{n+1}(t)|}{\delta^{k'}}.
$$
We estimate the second term in (\ref{prproof}) using (\ref{kappa}),
then (\ref{include}) and (\ref{goodPrime}) and finally
(\ref{snizu}), in the following way:
$$
|f_{k'}(c_{Q_t^{k'+1}})||(F^{[\varepsilon_{k'}]}_{Q_t^{k'+1}})'(t)|\leq
(1+\kappa)|f_{k'}(t)|\frac{c(\alpha)\lambda^\beta}{\delta^{k'+1}\tau^{\beta(n-k'+1)}}\leq$$
$$
\leq
(1+\kappa)(4/\theta)^{n-k'+1}|f_{n+1}(t)|\frac{c(\alpha)\lambda^\beta}{\delta^{k'+1}
\tau^{\beta(n-k'+1)}}.
$$
Finally, using (\ref{ghprime}) and then again (\ref{snizu}), we
get
$$
|T'_{l+1}(t)|\leq
\frac{c(\alpha)\lambda^\beta|f_l(t)|}{\rho\delta^{l+1}}
\leq(4/\theta)^{n+1-l}\frac{c(\alpha)\lambda^\beta|f_{n+1}(t)|}{\rho\delta^{l+1}}
$$
Taking $\delta<\theta\tau^\beta/(4A)$ and substituting all the estimates into (\ref{prproof}), we get
$$
|f'_{n+1}(t)|\leq
\frac{\lambda^\beta|f_{n+1}(t)|}{\delta^{n+1}}(c_1(\alpha)/A^{n-k'+1}+c(\alpha)
\frac{4(1+\kappa)}{\theta\tau^\beta}+\frac4{\rho\theta}\sum_{k=0}^{n-k'+1}A^{-k}).
$$
Taking $A>2$, we get (\ref{fprime}), if $c_1(\alpha)$ is
sufficiently large.

\bigskip\noindent{\bf8. The end of the proof.} In order to complete the
construction, we should define the sets $G_n$. Fix a parameter $\tau$
in order to satisfy inequality
$C(\tau)\sum_{n=1}^\infty\frac{1}{n^2}<\frac18$, where $C(\tau)$
is a constant from Lemma 5.

We define the sets $G_{n+1}$ as follows: $G_1:=H_1=\{I\}$; an
interval $Q\in H_{n+1}$ belongs to $G_{n+1}$, if $Q\subset V_n$,
\begin{equation}
M_Q(f_n)\leq K_n\eta^n,\label{G1}
\end{equation} where $K_n$ and $\eta$ will be defined later
and, besides,
\begin{equation}
\widetilde{D_n}(c_Q)\leq n/2.\label{G2}
\end{equation}
(in fact, of course, $\widetilde{D_n}$ is a constant on $Q$). Note
that the choice of $\tau$ and Lemma 5 guarantee that the total (for all
$n$'s) length of $Q\in H_{n+1}$, $Q\subset V_n$, not included in
$G_{n+1}$ because of violation of the condition (\ref{G2}), does not
exceed $\frac18$.

\begin{lemma}
There exists a constant $C'(\alpha)$ such that for all
sufficiently small $\lambda$ the following inequalities hold:
\begin{enumerate}\item If $Q\in G^g_{n+1}$, then
$$\int_Q|f_{n+1}|^p\leq X\int_Q |f_n|^p.$$
\item If $Q\in G^d_{n+1}$, then
$$\int_Q|f_{n+1}|^p\leq Y\int_Q |f_n|^p,$$
\end{enumerate}
where
$X:=\gamma(\lambda)^p(1+C'(\alpha)\lambda^\beta)^p(=(1+B\lambda)(1+C'(\alpha)\lambda^\beta)^p)$,
$Y:=(1+C'(\alpha)\lambda^\beta)^p$
\end{lemma}
{\bf Proof.} Letting
$P_{n+1}(t):=f_n(t)-f_n(c_{Q_t})e^{i\theta}F^{[\varepsilon]}_{Q_t}(t)$,
we get for $Q\in G_{n+1}^g$:
$$
\int_{Q}|f_{n+1}|^p = \int_{Q}|P_{n+1}+T_{n+1}|^p \leq
\int_{Q}|P_{n+1}|^p(1+\frac{|T_{n+1}|}{|P_{n+1}|})^p.
$$
Applying point 2 of Lemma 4 and the estimate (\ref{ghvost}),
we get
$$
\int_{Q}|f_{n+1}|^p\leq \int_{Q}
|P_{n+1}|^p(1+\frac{2c(\alpha)\lambda^\beta}{\rho\theta})^p,
$$
and then, estimating the integral using point 1 of Lemma 4, we
have
$$
\int_{Q}|f_{n+1}|^p \leq
(1+B\lambda)(1+C'(\alpha)\lambda^\beta)^p\int_{Q}|f_n|^p.
$$
Recall that applicability of Lemma 4 is provided by (\ref{kappa}).

The second case is even easier.
\begin{lemma}
If $\lambda$ is sufficiently small, then for all $n$ there holds
an inequality
\begin{equation} \int_{V_n}|f_n|^p\leq \eta^n\int_{V_1}
|f_1|^p\label{neqeta}
\end{equation}
with some $\eta\in (0,1)$.
\end{lemma}
{\bf Proof.} Denote as $\Theta_n$ the set $\{0,1\}^n$ (the set of
all ordered sets $v=(v_1,\dots,v_n)$ of the length $n$ of zeroes
and unities). Define for $v\in\Theta_n$ the set $G^{(v)}\subset
G_{n}$ as the set of all $Q$'s such that for all $k$ $v_k=1$ if
and only if $ Q\subset V_k^g$. Let for $v\in\Theta_n$
$Z(v):=\prod\limits_{k=1}^n X^{v_k} Y^{(1-v_k)}$. Let by
definition $G_1^g:=I$ and prove by induction in $n$ the estimate
\begin{equation} \int_I|f_1|^p\geq
X\sum_{v\in\Theta_n}Z(v)^{-1} \int_{G^{(v)}}|f_n|^p.
\label{derevo}
\end{equation}
The base is obvious. Suppose (\ref{derevo}) holds for some n. Then
$$
Z(v)^{-1} \int_{G^{(v)}}|f_n|^p=Z(v)^{-1}
(\int_{G^{(v,1)}}|f_n|^p+\int_{G^{(v,0)}}|f_n|^p) \geq
$$
$$
\geq
Z(v)^{-1}(X^{-1}\int_{G^{(v,1)}}|f_{n+1}|^p+Y^{-1}\int_{G^{(v,0)}}|f_{n+1}|^p)=
$$
$$
=Z(v,1)^{-1}\int_{G^{(v,1)}}|f_{n+1}|^p+Z(v,0)^{-1}\int_{G^{(v,0)}}|f_{n+1}|^p.
$$
(we used Lemma 7). Thus we have proved the inductive step.
Then, for $\lambda$ sufficiently small we have $X<1$, $Y>1$. If
$n>N_0$, $v\in \Theta_n$, then $\sum_k v_k < n/2$ implies
$G^{(v)}=\varnothing$ (it follows from (\ref{G2})), therefore the
right-hand side in (\ref{derevo}) is not less than
$X^{-\frac{n}2+1}Y^{-\frac{n}2}\int_{V_n}|f_n|^p$, and
($\ref{neqeta}$) follows provided $\lambda$ is sufficiently small.

Now we are ready to finish the setup of the construction. What we
should do is to make precise the condition (\ref{G1}). We take
$\eta$ from Lemma 8. As, by Lemma 8, $\int_{V_n}|f_n|^p\leq
C\eta^n$, the total length of all intervals $Q\in H_{n+1}$, such
that $Q\subset V_n$ and on $Q$ the condition (\ref{G1}) fails,
does not exceed $\frac{C}{K_n}$. For $K_n$ we take a sequence
growing as some power of $n$ and satisfying the condition
$\sum_{n=1}^\infty \frac{C}{K_n}<\frac18$. We get that total
length of all intervals dropped out of $V_n$ on all steps because
of violation of (\ref{G1}) is less than $\frac18$. But one can
tell the same thing about the length of all intervals dropped out
because of violation of (\ref{G2}). Hence $|\bigcap_{i=1}^\infty
V_n|>\frac34$.

It follows from (\ref{G1}) and an estimate (\ref{kappa}) of the
oscillation of $f_n$ that if $Q\in G_{n+1}$,  then $|f_n(c_Q)|\leq
C'\eta'^n$. Hence, taking into account (\ref{estimF}) and
(\ref{sigma}) we get
\begin{equation}\label{prirost} |f_n(t)-f_{n+1}(t)|
\leq
\frac{C_1\eta'^n}{\varepsilon_n^\beta}+C_2\lambda^\beta\eta'^n,
t\in\mathbb{R}
\end{equation}
The first term in the right-hand side corresponds to the building
block $F^{[\varepsilon_n]}_{Q^{n+1}_t}$ (if there is any), the
second -- to all the others. It follows from this estimate that,
$f_n$ converges uniformly on ${\mathbb R}$ to some function $f$.
We also need an estimate
\begin{equation}
\label{major} |f_n(t)-f_{n+1}(t)|\leq
c\frac{\lambda^\beta\delta_{n+1}^\alpha
C'\eta'^n}{|t|^\beta},\qquad t\notin 3I,
\end{equation}
which follows from (\ref{Fdaleko}). Now, (\ref{major}) implies
that
\begin{equation}
\label{maj} |f_n(t)|\leq c|t|^{-\beta},\quad |t|\geq 3/2,
\end{equation}
so we can, fixing $t$, write
$$
\int\limits_{\mathbb R}
\frac{f_n(s)ds}{|s-t|^{1-\alpha}}=\int\limits_{|s|\leq
\max(2|t|,3/2)}\frac{f_n(s)ds}{|s-t|^{1-\alpha}}+\int\limits_{|s|\geq
\max(2|t|,3/2)}\frac{f_n(s)ds}{|s-t|^{1-\alpha}}.
$$
In the first term, the passage to the limit in the integral is
provided by uniform convergence of $f_n$ as  $n\rightarrow\infty$,
in the second one, the function under integral is majorized by
$c|s|^{-2}$. Hence
$$
\lim\limits_{n\rightarrow\infty}g_n =
\lim\limits_{n\rightarrow\infty}\int\limits_{\mathbb R}
\frac{f_n(s)}{|s-t|^{1-\alpha}}=\int\limits_{\mathbb R}
\frac{f(s)}{|s-t|^{1-\alpha}}
$$
(The first inequality follows from Lemma 1). Finally note that
the functions $g_n$ converge uniformly to $g$:
$$
|r_n|\leq \frac{C''\eta'^n}{\varepsilon_n}.
$$
So, the program announced in the beginning of Section 3 is
completed. Let us now show that the function $f$
satisfies H\"{o}lder's condition.

In this argument, constants may depend on all parameters except
$n$. Of course, we use the proposition from the beginning of Section
4. The condition $|f_n-f|\leq C_1\eta_1^n$ follows easily from
(\ref{prirost}).

Let us estimate the derivative of the correcting term on the
$n$-th step:

$$
|(W_\alpha r_n)'(t)|=|\frac{d}{dt}(\sum\limits_{Q\in
G^g_{n+1}}f_n(c_Q)F^{[\varepsilon_n]}_Q(t)e^{i\theta})|\leq
\max\limits_{Q\in G^g_{n+1}}|f_n(c_Q)|\sum\limits_{Q\in
G^g_{n+1}}|(F^{[\varepsilon_n]}_Q(t))'|\leq
$$
$$
\leq\max\limits_{Q\in
G^g_{n+1}}|f_n(c_Q)|\left(|(F^{[\varepsilon_n]}_{Q_t^{n+1}}(t))'|+\sum\limits_{Q\in
H_{n+1},\:Q\neq Q_t^{n+1}}|(F^{[\varepsilon_n]}_Q(t))'|\right).
$$

We know that $\max\limits_{Q\in G^g_{n+1}}|f_n(c_Q)|\leq
C'\eta'^n$. Then, we use
(\ref{estimFP}) for the first term in the brackets, and for the sum -- (\ref{sigmaP}) with
$\epsilon:=\frac{\delta^{n+1}}2$. Thus we have
$$|(W_\alpha
r_n)'(t)|\leq C'\eta'^n(\frac
{C_1}{\delta^{n+1}\varepsilon_n}+\frac{C_2}{\delta^{n+1}})\leq
\frac {C_3}{\delta'^{n+1}}$$

Hence $|f_n'(t)|\leq
\frac{C_3}{1-\delta'}\delta'^{-n-1}\leq\frac{C_4}{\delta'^{n+1}}$,
and H\"{o}lder's property for $f$ is proved.

\bigskip\noindent{\bf9. Remark on the order of choice of the parameters.}
Recall the order we chose our parameters in. Given $\alpha\in
(0,1)$, we fix $p$ and $\theta$ (Lemmas 3 and 4). Here Lemma 4
holds for all sufficiently small $\lambda$ and $\varepsilon$. Then
we fix a sequence $\varepsilon_n$. Independently of other
parameters we fix $\tau$. Lemma 6 holds for all sufficiently small
$\lambda$ and $\delta$, independently of the choice of $G_n$ (here how
small $\lambda$ and $\delta$ should depend on $\rho$ from this
lemma, and $\rho$ itself only depends on $\alpha$, $\theta$, $p$,
$\tau$, and the setup function $g_1$). Now we choose $\lambda$ such
that the multiplier $\eta$ in front of the integral in the
right-hand side of (\ref{neqeta}) is less than 1 (decreasing
$\lambda$ killing the influence of ``tails'' $T_n$ in comparison
with the correction effect provided by point 1 of Lemma 2. Note
that the last effect decays when $\lambda$ decreases as well
($\gamma(\lambda) \stackrel{\lambda\rightarrow
0}{\longrightarrow}1$!), but the ``tails'' dies faster). Finally,
we fix $\delta$ in order to provide (\ref{kappa}).

So far our considerations did not depend on $G_n$, therefore
we did not need to define them. Now we fix a sequence $K_n$ and
this finishes the definition of our construction.

\bigskip\noindent{\bf10. The case of negative $\alpha$.} One can
ask whether there is an analogue of Theorem 1 for other
M.Rietz's kernels. In the case of the kernel
$|x|^{-\beta},\quad1<\beta<2$, the answer is affirmative; moreover,
as the convolution with such a kernel is, at least formally, the
inverse operator to $U_{\beta-1}$ (see Lemma 1, points 2 and 5),
the example in essential coincides with the one built above.
\begin{Theorem}
There exist a nonzero continuous function $g:{\mathbb
R}\rightarrow {\mathbb C},$ $supp\:g\subseteq 3I$ and a set $E$ of
positive measure, such that $t\in E \Rightarrow g(t)=0,$
$\int_{\mathbb R} g(x)|t-x|^{-\beta}dx=0$. The last integral
converges absolutely for every $t\in E$. The function $g$
satisfies H\"{o}lder's condition with an exponent $\beta-1$.
\end{Theorem}
{\bf Proof.} We take for $g$ the function built in Theorem 1 (with
$\beta-1$ for $\alpha$). Let
$$
E:=V\cap\left(I\backslash S\right),\quad\mbox{where
}S:=\bigcup\limits_{n\in {\mathbb N}}\bigcup\limits_{Q\in
G^g_{n+1}}3\varepsilon_n(Q-c_Q)+c_Q
$$
The idea is that $|S|< \frac34$, but now $supp\; g_n$ is contained
``rather deep'' inside $S$:
\begin{equation}
\label{supp} dist(supp\;g_n,E)\geq\varepsilon_n\delta^{n+1}.
\end{equation}

We know (point 5 of Lemma 1), that for $t\in E$
$$
f_n(t)=\int_{\mathbb R}
g_n(x)|t-x|^{-\beta}dx\stackrel{n\rightarrow
\infty}{\longrightarrow}0.
$$
To prove Theorem 2, it is sufficient to justify the passage to
the limit in the integral. For the summable majorant we take a
function
$$
\widetilde{g}(x):=|t-x|^{-\beta}\sum\limits_{n=1}^\infty|g_{n+1}(x)-g_n(x)|.
$$
Let us estimate the $n$-th term:
$$
|g_{n+1}-g_n|\leq\sum_{Q\in
G^g_{n+1}}|f_n(c_Q)|(\delta^{n+1}\lambda)^\alpha(\phi_{\varepsilon_n
Q})\leq C'\eta'^n (\delta^{n+1}\lambda)^\alpha\sum_{Q\in
G^g_{n+1}}(\phi_{\varepsilon_n Q}).
$$
Check summability of the majorant:
$$
\int_{\mathbb R} |t-x|^{-\beta}|g_{n+1}(x)-g_n(x)|\leq C'\eta'^n
(\delta^{n+1}\lambda)^\alpha\sum_{Q\in G^g_n}\int_{\mathbb
R}|t-x|^{-\beta}\phi_{\varepsilon_n Q}(x)dx\leq
$$
$$ \leq C'\eta'^n (\delta^{n+1}\lambda)^\alpha\sum_{Q\in
G^g_n}|t-x^*_Q|^{-\beta}\int_{\mathbb R}\phi_{\varepsilon_n
Q}(x)dx\leq
$$
\begin{equation}\leq C'\eta'^n
(\delta^{n+1}\lambda)^\beta\sum_{Q\in G^g_n}|t-x^*_Q|^{-\beta}.
\label{ocmajor}
\end{equation}
Here $x^*_Q$ denotes a point of the support of
$\phi_{\varepsilon_n Q}(x)$, closest to $t$. Using (\ref{supp}) and
the fact that the distance between two different points $x^*_Q$ and
$x^*_{Q'}$ is no less than $\delta^{n+1}$, we get
$$
\sum\limits_{Q\in G^g_n}|t-x^*_Q|^{-\beta}\leq
2\sum\limits_{k=0}^\infty(\varepsilon_n\delta^{n+1}+k\delta^{n+1})^{-\beta}\leq
2\delta^{(n+1)(-\beta)}\sum\limits_{k=0}^\infty|\varepsilon_n+k|^{-\beta}\leq
$$
\begin{equation}
\leq
2\delta^{(n+1)(-\beta)}(\varepsilon_n^{-\beta}+\sum\limits_{k=1}^\infty
k^{-\beta})\leq C\delta^{(n+1)(-\beta)}\varepsilon_n^{-\beta}
\label{ocsum}
\end{equation}
Substituting (\ref{ocsum}) to (\ref{ocmajor}) and summing over all
$n$, we get
$$
\int_{\mathbb R}\widetilde{g}(t)dt\leq
C''\lambda^{-\beta}\sum\limits_{n=1}^\infty\eta'^n\varepsilon_n^{-\beta}<+\infty
$$
We should now only prove that $g$ satisfies H\"{o}lder's condition
with an exponent $\alpha=\beta-1$. In fact this is a property of
the potential $U_{\alpha}$ of any bounded function for which it is
defined. To prove it, take $t>0$ and write the following estimate:
$$
\int_{\mathbb R}||x|^{\alpha-1}-|x-t|^{\alpha-1}|dx =
\int_{(-t;2t)}+\int_{{\mathbb R}\backslash(-t;2t)}=:J_1+J_2.
$$
Estimate each term:
$$
J_1 \leq \int_{(-t;2t)}|x|^{\alpha-1} +
\int_{(-t;2t)}|x-t|^{\alpha-1} = \frac2\alpha(1+2^\alpha)t^\alpha;
$$
$$
J_2 \leq 2\int_{(t;+\infty)}x^{\alpha-1} -(x+t)^{\alpha-1} \leq
2(\alpha-1)\int_{(t;+\infty)}tx^{\alpha-2} \leq 2t^\alpha.
$$
When estimating $J_2$ we have used the inequality
$|h(x+t)-h(x)|\leq t\sup\limits_{s\in (x,x+t)}|h'(s)|$ for a
smooth function $h$. From this estimates we get
$$
|(U_\alpha f)(t+\delta)-(U_\alpha f)(t)| \leq \sup\limits_{\mathbb
R}|f|\int_{\mathbb
R}(|t+\delta-x|^{\alpha-1}-||t-x|^{\alpha-1}|)dx \leq
C(\alpha)(\sup\limits_{\mathbb R}|f|)\delta^\alpha.
$$
The theorem is proven.

{\bf Remark 1.} Proving the smoothness of $g$ we didn't use all
the information we had about $f$. In fact $f$, besides it is
bounded and belongs to the domain of $U_\alpha$, satisfies
H\"{o}lder's condition with some exponent $r>0$. Using this and
well-known techniques of estimating operators similar to M. Rietz
potential (see, for example, \cite{Zyg}), one can
prove that $g$ satisfies H\"{o}lder's condition with an exponent
$\beta-1+r$.

{\bf Remark 2.} The theorem proven in \cite{Havin} states that for
the potentials $U_\alpha$, $-1<\alpha<0$, uniqueness holds if the
density $g$ belongs to $C^{1+\varepsilon}$ with some
$\varepsilon>0$. Theorem 2 shows that the last condition cannot be
replaced by H\"{o}lder's condition with an exponent $-\alpha$. So,
there is a gap between the two results, which decreases when
$\alpha\rightarrow -1$. If $\alpha\leq-1$, one does not need any
supplementary smoothness condition: the mere existence of the
potential is sufficient (see, for example, \cite{HJ}).

For the cases $\alpha=0$ and $\alpha>1$ (except odd integers for
which the uniqueness does not hold in any sense), the question
whether one can omit the smoothness conditions imposed
in \cite{Havin} remains open.

\bigskip\noindent{\bf11. Extension of the results to the multidimensional case.}
The result of Theorem 1 can be extended to the case of M. Rietz
potentials in spaces $\mathbb{R}^d$ for $d>1$. In this case
for $\alpha\in(0,d)$  we consider a set of all measurable
functions $f$, satisfying the condition
\begin{equation}
\label{Cond1Rn}
\int_{\mathbb{R}^d}\frac{|f(x)|}{1+|x|^{d-\alpha}}dx < + \infty
\end{equation}
(as above, we denote this set $dom\,U_\alpha$). We let
$$
U_\alpha f := f\ast |x|^{d-\alpha},\; f\in dom\, U_\alpha,
$$
where $\ast$ denotes the convolution in $\mathbb{R}^d$. The case
of major interest is $d=2$; $\alpha=1$ (Newton's potential of the
charge concentrated in the plane). Note, however, that in case of
$d>1$ there is not any analogue of the uniqueness theorem
mentioned in the introduction.

There holds the following generalization of Theorem 1:

\begin{Theorem}
For all $d\in\mathbb{N}$ and for all $\alpha\in(0,d)$ there exist a
nonzero function $f:\mathbb{R}^d\rightarrow \mathbb{R}$, $f\in
dom\,U_\alpha$ and a set $E\subset \mathbb{R}^n$ of positive
Lebesgue measure satisfying the condition $f|_E=0$, $U_\alpha
f|_E=0$, and H\"{o}lder's condition with some positive exponent.
\end{Theorem}

The proof of this theorem is quite similar to the one of Theorem 1. We
highlight some details differing in the multidimensional case.

We need an operator $W_\alpha$, ``the inverse operator'' to
$U_\alpha$. The precise expression of this operator (see, for
example, \cite[page 241]{BH}) does not matter for us, the only
thing we need is that if we now denote as $\beta$ the number
$d+\alpha$, then points 1,2,3,5 of Lemma 1 still hold.

The role of $I$ will be played by the cube $I^d$, and we shall
consequently divide it to congruent cubes with the side equal to
$\delta^n$. Instead of ``the finitizator'' $\phi(x)$ we take the
function $\phi(|x|)$

The computations made in \cite[page 241]{BH} show that lemma 2
still holds in the multidimensional case. Lemmas 3 and 4 can be derived
from it quite similarly to the above.

Taking into account that now $\beta=d+\alpha$, the most of
computations in the multidimensional case will repeat one-dimensional
literally, if we also replace derivative by gradient
everywhere. So, because of point 3 for lemma 1 there still hold (\ref{estimF}), (\ref{Fdaleko}), (\ref{FPdaleko}),
(\ref{goodF}), (\ref{goodPrime}) (now for the cube $Q$ with sides
parallel to the coordinate axes the symbol $Q[a]$ denotes $Q\backslash
Q'$, where $Q'$ is a cube obtained from $Q$ by homothety with the
center $c_Q$ and the dilation factor $a$).

The remaining part of the construction is the same. Note
that estimates (\ref{sigma}) and (\ref{sigmaP}), playing the key
role in the proof of lemma 6, and seeming to depend
on the dimension, in fact hold in the literally same form.

{\bf Remark.} Theorem 2 can be extended to the multidimensional
case as well: \textit{ for $d<\beta<2d$ there exist a nonzero
continuous function $g:{\mathbb{R}^d}\rightarrow {\mathbb C}$ and a
set $E$ of positive measure, such that if $t\in E$, then $g(t)=0$,
$\int_{\mathbb{R}^d} g(x)|t-x|^{-\beta}dx=0$. The last integral
converges absolutely for all $t\in E$. The function $g$ satisfies
H\"{o}lder's condition with an exponent
$\min\{\beta-d;1\}$}\footnote{This smoothness estimate can be
obtained by a simple method similar to the one used in the proof of
Theorem 2. Using the techniques mentioned in the remark after
Theorem 2, one can prove the inclusion $g\in C^{\beta-d+r}$ (in
case when $\beta-d+r$ is an integer, we understand it as
the corresponding Zygmund class)}. The only difference in the proof is
the estimate (\ref{ocsum}), where sums becomes multiple (of order
$d$).  They still will converge, because now $\beta>d$.



\begin{thebibliography}{99}
\bibitem{AK}
Aleksandrov, A. A., Kargaev, P. P. Hardy classes of functions harmonic in half-space, Algebra i Analiz 5:2 (1993), page 1-73(Russian).
\bibitem{BH}
Belyaev, D. B.  and Havin, V. P. On the uncertainty principle for
M.Rietz potentials., Arkiv f\"{o}r Mathematik, 39(2001), 229-233
\bibitem{Binder}
Binder, I., Theorem on correction by harmonic gradients, Algebra i Analiz 5:2 (1993), page 91-107(Russian).
\bibitem{BW}
Bourgain, J., Wolff, T., A remark on gradients of harmonic
functions in dimention $\geq3$, Colloq. Math. 40/41 (1990), no. 1,
25-260.
\bibitem{Havin}
Havin, V. P., Uncertainty principle for one-dimensional M. Rietz potentials
Dokl. AN SSSR 264, $N^o$3 (1982), 559-563.(Russian)
\bibitem{HJ}
Havin, V., J\"{o}ricke, B., The Uncertainty Principle in Harmonic
Analysis, Springer-Verlag, Berlin,1994.
\bibitem{wolff}
Wolff, T., Counterexamples with harmonic gradients in
$\mathbb{R}^3$, Essays on Fourier Analysis in Honor of Elias M.
Stein (Princeton, N. J., 1991) (Fefferman, C., Fefferman, R., and
Wainger, S., eds.) 321-384, Princeton Univ. Press, Princeton,
1995.
\bibitem{Zyg}
Zygmund A., Trigonometric Series, Cambridge University press, 2002.
\end{thebibliography}
\end{document}